\documentclass[11pt,a4paper]{article}
\usepackage[english]{babel}
\usepackage{latexsym}
\def\qed{\vbox{\hrule
\hbox{\vrule\hbox to 5pt{\vbox to 8pt{\vfil}\hfil}\vrule}\hrule}}

\newcommand{\beg}{\begin{eqnarray*}}
\newcommand{\begn}{\begin{eqnarray}}
\newcommand{\en}{\end{eqnarray*}}
\newcommand{\enn}{\end{eqnarray}}

\begin{document}
\vspace*{0.5cm} \begin{center} \noindent {\LARGE \bf Compactness theorems for gradient Ricci solitons }\\[0.7cm]
Xi Zhang\footnote[1]{The author was supported by NSF in China, No.10201028}\\[0.2cm]
Department of Mathematics, Zhejiang University,\\
 Hangzhou 310027, Zhejiang,
P. R. China \\
E-mail: xizhang@zju.edu.cn\\[1cm]\end{center}

{\abstract In this paper, we prove the compactness theorem for
gradient Ricci solitons. Let $(M_{\alpha} , g_{\alpha})$ be a
sequence of compact gradient Ricci solitons of dimension $n\geq
4$, whose curvatures have uniformly bounded $L^{\frac{n}{2}}$
norms, whose Ricci curvatures are uniformly bounded from below
with uniformly lower bounded volume and with uniformly upper
bounded diameter, then there must exists a subsequence
$(M_{\alpha}, g_{\alpha})$ converging to a compact orbifold
$(M_{\infty}, g_{\infty})$ with finitly many isolated
singularities, where $g_{\infty}$ is a gradient Ricci soliton
metric in an orbifold sense.
\\[0.5cm]
{\bf Subj.Class:} differential geometry.

\section{Introduction}
\setcounter{equation}{0}

\hspace{0.4cm}

The concept of convergence of Riemannian manifolds was introduced
by M.Gromov ([Gr]). As is now well-known, the Cheeger-Gromov
convergence theorem ([Ch], [Gr], [GW], [Ka], [Pe]) implies that
the space $\mu (\Lambda , \nu , D)$ of compact Riemannian
$n$-manifolds of sectional curvature $|K|\leq \Lambda $, volume
$\geq \nu >0$ and diameter $\leq D$, is precompact in the $C^{1,
\alpha }$ topology. There has been increasing interest lately in
compactness theorems of Riemannian manifolds under various
geometric assumptions ([An1], [An2], [Ga], [Na], [Ti1], [Ti2]).
For instance, in [An1] and [Na], the authors show that if $\{
(M_{\alpha}, g_{\alpha })\}$ is a sequence of Einstein manifolds
of dimension $n$ satisfying: $diam (M_{\alpha}, g_{\alpha})\leq
C$; $\int_{M_{\alpha}}
\|Rm(g_{\alpha})\|_{g_{\alpha}}^{\frac{n}{2}}\,
dV_{g_{\alpha}}\leq C$; and $Vol(M_{\alpha}, g_{\alpha})\geq
\frac{1}{C}$, where $C$ is a uniform constant, then there is a
subsequence of $\{ (M_{\alpha}, g_{\alpha})\}$ converges to an
Einstein orbifold with finitely many isolated singular points.
Also see [Ti1] and [Ti2] for the case of K\"ahler-Einstein
manifolds.

The Ricci flow equation
\begin{eqnarray}
\frac{d}{dt} g_{ij}=-2R_{ij}
\end{eqnarray}
has been introduced by R.Hamilton in his seminal paper [Ha1].
Natural questions that arises in studying the Ricci flow equation,
are under which conditions a solution will exist for all times, if
there exists a limit of the solutions when we approach infinity
and how we can describe the metrics obtained in the limit. In the
case of dimension three with positive Ricci curvature and
dimension four with positive curvature operator we know (due to
R.Hamilton [Ha1], [Ha2]) that the solutions of the Riccci flow ,
in both cases exist for all times, converging to Einstein metrics.
In general, we can not expect to get an Einstein metric in the
limit, but we can expect to get solitons in the limit. In brief, a
soliton is just a solution to the Ricci flow (1.1)  which moves by
diffeomorphisms and also shrinks or expands by a factor at the
same time. Such a solution is called a homothetic shrinking or
expanding Ricci soliton. The equation for a homothetic Ricci
soliton is
\begin{eqnarray}
2\lambda g_{ij}
-2R_{ij}=g_{ik}\nabla_{j}v^{k}+g_{jk}\nabla_{i}v^{k},
\end{eqnarray}
where $\lambda $ is the homothetic constant, $v$ is the vector
field induced by the $1$-parameter family of diffeomorphisms. For
$\lambda >0$ the soliton is shrinking, for $\lambda <0$ it is
expanding, and the case $\lambda =0$ is a steady Ricci soliton. If
the vector field $v$ is the gradient of a function $u$ we say that
the soliton is a gradient Ricci soliton, thus
\begin{eqnarray}
\lambda g_{ij} -R_{ij}=\nabla_{i}\nabla_{j} u,
\end{eqnarray}
is the gradient Ricci soliton equation. Einstein metrics  can be
considered as Ricci soliton when the vector field $v$ is zero.

In this paper, we want to consider the compactness result for the
Ricci soliton. When the underly manifolds are closed (compact,
without boundary), one can easily check that the steady and
expanding Ricci solitons are in factly Einstein metrics. So, we
mainly consider the shrinking case. We prove the following
theorem.

\hspace{0.3cm}

{\bf Theorem 1.1} {\it Let $(M_{\alpha}, g_{\alpha})$ be a
sequence of shrinking gradient Ricci solitons of dimension $n\geq
4$, i.e, satisfy the following equation
\begin{eqnarray}
g_{\alpha} -Ric(g_{\alpha})=\nabla d u_{\alpha},
\end{eqnarray}
such that

(1),  $Ric(g_{\alpha})\geq -C_{1}g_{\alpha}$;

(2), $diam (M_{\alpha}, g_{\alpha})\leq C_{2}$;

(3), $Vol (M_{\alpha }, g_{\alpha })\geq C_{3}$;

(4), $\int_{M_{\alpha}}|Rm|^{\frac{n}{2}} \, dV_{g_{\alpha}}\leq
C_{4}$;

for some uniform positive constants $C_{1}$, $C_{2}$, $C_{3}$,
$C_{4}$. Then there is a subsequence $(M_{\alpha}, g_{\alpha})$
converging to $(M_{\infty}, g_{\infty})$ in Cheeger-Gromov's
sense, where $M_{\infty}$ is an orbifold with finitely many
isolated singularities and $g_{\infty}$ is a Ricci soliton in an
orbifod sense.

Further, if $n$ is odd, there are no singular points and
$(M_{\infty}, g_{\infty})$ is a smooth gradient Ricci soliton
which diffeomorphic to $M_{\alpha}$, for $\alpha $ sufficiently
large. In this case, $(M_{\alpha}, g_{\alpha})$ (sub)converges
smoothly to $(M_{\infty}, g_{\infty})$.}

\hspace{0.3cm}

{\bf Remark 1.2:} {\it An $n$-dimensional orbifold $M_{\infty}$ is
a topological space satisfying: each point $x$ in $M_{\infty}$
admits an open neighborhood $U_{x}$ homeomorphic to
$B^{n}/\Gamma_{x}$, where $B^{n}$is the unit disc in $R^{n}$ and
$\Gamma_{x}\subset O(n)$ is a finite group and those $U_{x}$ are
patched together by smooth transition functions. Any point $x$
with $\Gamma_{x}$ trivial is called a regular point of
$M_{\infty}$. In particular, $M_{\infty}$ is a manifold near such
a regular point. Denote by $Reg(M_{\infty})$ the set of all
regular points. All other points are singular points of
$M_{\infty}$ i.e. $Sing(M_{\infty})=M_{\infty}\setminus
Reg(M_{\infty})$. We will confine ourselves to the special case
that $Sing(M_{\infty})$ consists of isolated points. A Ricci
soliton $g_{\infty}$ in an orbifold sense is just the one on
$Reg(M_{\infty})$ such that for each $x\in Sing(M_{\infty})$ if
$\pi_{x}: B^{n} \rightarrow U_{x}$ is the local uniformization,
then there is a diffeomorphism $\varphi$ of $B^{n}$ such that
$\varphi^{\ast}\pi_{x}^{\ast}g_{\infty}$ can be extended smoothly
to a gradient Ricci soliton $C^{\infty}$-metric on $B^{n}$.

In theorem 1.1, we say that $(M_{\alpha}, g_{\alpha})$ converge to
an orbifold $(M_{\infty}, g_{\infty})$ in Cheeger-Gromove's sense,
if for any compact subset $K\subset M_{\infty}\setminus
Sing(M_{\infty})$ there are compact sets $K_{\alpha}\subset
M_{\alpha }$ and diffeomorphisms $\phi_{\alpha} : K\rightarrow
K_{\alpha}$ so that $\phi_{\alpha}^{\ast}g_{\alpha}$ converge to
$g_{\infty}$ in $C^{\infty}$ topology.}

\hspace{0.3cm}

Recall the formula of Avez ([Be]) expressing the Euler
characteristic $\chi (M)$ of a compact $4$-manifold in terms of a
curvature integral
\begin{eqnarray}
\chi (M) =\frac{1}{8\pi ^{2}}\int_{M} |Rm|^{2}-4|Ric |^{2}+R^{2} ,
\end{eqnarray}
where $R$ is the scalar curvature. Clearly, a bound on $\int_{M}
|Ric|^{2}$ and the second Betti number $b_{2}(M)$ implies a bound
on $\int_{M} |Rm|^{2}$. The Bishop comparison theorem implies
there is a upper bound of volume by the lower bound of Ricci
curvature and the upper bound of diameter. On the other hand, in
section two, we will prove that when $g$ is a gradient Ricci
soliton, the lower bound of Ricci curvature, the upper bound of
diameter, and the lower bound of volume imply an upper bound of
scalar curvature, then we have an upper bound of Ricci curvature.
So, for Ricci solitons, a lower bound of Ricci curvature and
volume, an upper bound of diameter, and a  bound of $b_{2}(M)$
imply a bound on $\int_{M} |Rm|^{2}$. We have the following
corollary.

\hspace{0.3cm}

{\bf Corollary 1.3} {\it Let $(M_{\alpha}, g_{\alpha})$ be a
sequence of shrinking gradient Ricci solitons of dimension $ 4$,
such that

(1),  $Ric(g_{\alpha})\geq -C_{1}$;

(2), $diam (M_{\alpha}, g_{\alpha})\leq C_{2}$;

(3), $Vol (M_{\alpha }, g_{\alpha })\geq C_{3}$;

(4), $b_{2}(M_{\alpha})\leq C_{4}$,

for some uniform positive constants $C_{1}$, $C_{2}$, $C_{3}$,
$C_{4}$. Then there is a subsequence $(M_{\alpha}, g_{\alpha})$
converging to $(M_{\infty}, g_{\infty})$ in Cheeger-Gromov's
sense, where $M_{\infty}$ is an orbifold with finitely many
isolated singularities and $g_{\infty}$ is a Ricci soliton in an
orbifod sense.}

\hspace{0.3cm}

More recently, H.D.Cao and N.Sesum ([CN]) proved a compactness
result for the K\"ahler Rici solitons, where  the upper bound of
diameter was replaced  by a uniform lower bound of Perelman's
functional $\mu (g , \frac{1}{2})$ ([Pe1]). They point out that by
the same proof as in [Pe2] they can show that the uniform lower
bound of Ricci curvature, lower bound of Perelman's functional
$\mu (g , \frac{1}{2})$ , and the Euclidean volume growth implies
a uniform bound of diameter. In section two, we can show that, for
a sequence of gradient Ricci soliton, an uniform lower bound of
Ricci curvature and volume, and an uniform bound of diameter will
give an uniform bound of Perelman's functional $\mu (g ,
\frac{1}{2})$.

The key analytic tool in this paper is Uhlenbeck's ([Uh])
Yang-Mills estimate for curvatures of Yang-Mills connections. The
method to prove theorem 1.1 is similar with that in Anderson's
paper ([An1]) for Einstein case. The main point is obtaining the
$\epsilon$-regularity estimates for Ricci soliton which say that a
smallness of the $L^{\frac{n}{2}}$ norm of curvature implies a
point-wise bound on the curvature. Moreover, different with the
Einstein case, we should to obtain $C^{1}$ bound for potential
functions $u_{\alpha}$. It is well know ([Ya]), a lower bound for
the Ricci curvature and volume and an upper bound on the diameter
give a lower bound for the Sobolev constant $C_{s}$ of compact
manifold $M$,
\begin{eqnarray}
\|u\|_{\frac{2n}{n-2}}\leq
\frac{1}{C_{s}}\|du\|_{2}+Vol(M)^{-\frac{2}{n}}\|u\|_{2},
\end{eqnarray}
for any Lipschitz function $u$ on $M$. By the above Sobolev
inequality, we can  use the Moser iteration argument to obtain the
above estimates.

The organization of this paper is as follows. In section two, we
deduce some estimates, specially, we obtain the $C^{1}$ estimates
for functions $u_{\alpha}$, and an uniform bound for Ricci
curvature and Perelman's function. In the section three, we obtain
the $\epsilon $-regularity estimate for Ricci soliton. In the
section four and five, we give the proof of theorem 1.1.

\hspace{0.3cm}

\section{Preliminary Results}
\setcounter{equation}{0}

 \hspace{0.3cm}

Let $M$ be a compact manifold without boundary, and $g$ be a
gradient Ricci soliton i.e. satisfies formula (1.3), here we
assume that $u$ satisfying
\begin{eqnarray}
(2\pi )^{-\frac{n}{2}}\int_{M} e^{-u} \, dV_{g} =1
\end{eqnarray}
and  $Ric(g)\geq -C_{1}g$; $diam (M, g)\leq C_{2}$; $Vol (M,
g)\geq C_{3}>0$.

From formula (2.1), we have
\begin{eqnarray}
\inf_{x\in M} u \leq ln Vol(M,g) -\frac{n}{2} ln (2\pi).
\end{eqnarray}
On the other hand, the Bishop comparison theorem implies there is
an upper bound of volume by a lower bound of Ricci curvature and
an upper bound of diameter. So, there exists a constant $C_{5}$
depending only on $C_{1}$ and $C_{2}$ such that
\begin{eqnarray}
\inf_{x\in M} u \leq C_{5}.
\end{eqnarray}

Let $f=e^{-\frac{u}{2}}$, then we have
\begin{eqnarray}
\begin{array}{lll}
\triangle f^{2} &=e^{-u}|\nabla u|^{2}-e^{-u}\triangle u\\
&=e^{-u}|\nabla u|^{2}+e^{-u}(R-n\lambda )\\
&\geq 4|\nabla f|^{2}-nf^{2}(C_{1}+\lambda).
 \end{array} \end{eqnarray}
From the above Bochner type inequality and the Sobolev inequality
(1.6), using the Moser iteration argument ([Zh], proposition2.2),
we have the following mean value inequality,
\begin{eqnarray}
\sup_{x\in M} f \leq C_{6} (\int_{M} e^{-u} dV_{g})^{\frac{1}{2}},
\end{eqnarray}
where $C_{6}$ depend only on $C_{1}$, $C_{2}$ and $C_{3}$. So, we
obtain a lower bound of $u$, i.e. there is a constant $C_{7}$
depending only on $C_{1}$, $C_{2}$ and $C_{3}$ such that
\begin{eqnarray}
\inf_{M} u \geq -C_{7}. \end{eqnarray}

From the Ricci soliton equation and the lower bound of Ricci
curvature, we have
\begin{eqnarray} \nabla d u \leq
(C_{1}+\lambda)g.
\end{eqnarray}
Let $P, Q\in M$ such that $u(P)=\inf_{x\in M}u, u(Q)=\sup_{x\in
M}u$, and $\gamma : [0, d]\rightarrow M$ be a minimizing geodesic
connecting $P$ and $Q$, i.e. $\gamma (0)=P$, $\gamma (d) =Q$, here
$d=dist (P, Q)$. We have
\begin{eqnarray}
\begin{array}{lll}
\frac{d u(\gamma (t))}{d t}&=\langle \nabla u , \gamma '  \rangle_{\gamma (t)}-\langle \nabla u , \gamma ' \rangle_{p}\\
&=\int_{0}^{t}\frac{\partial }{\partial s}(\langle \nabla u , \gamma ' \rangle_{\gamma (s)})\, ds\\
&=\int_{0}^{t}\nabla_{\gamma '(s)}\langle \nabla u , \gamma '
\rangle \, ds \\
&=\int_{0}^{t}(\nabla du ) ( \gamma ', \gamma ' ) \, ds\\
&\leq (C_{1}+\lambda )t \\,
 \end{array} \nonumber\end{eqnarray}
and
\begin{eqnarray}
u(Q)-u(P)=\int_{0}^{d}\frac{d u(\gamma (t))}{d t} \, dt\leq
\int_{0}^{d} (C_{1}+\lambda )t  dt =\frac{1}{2}(C_{1}+\lambda
)d^{2}. \nonumber\end{eqnarray}

 From the above inequality and (2.3), we know that there exist a
 constant $C_{8}$ depending only on $C_{1}$, $C_{2}$, and $C_{3}$, such that
\begin{eqnarray}
\sup _{x\in M} u \leq C_{8}. \end{eqnarray}

Next, we want to obtain the estimate of $|\nabla u|$. From
equation (1.3), we have
\begin{eqnarray}
\nabla_{i} R_{jk} -\nabla_{j} R_{ik} =-R_{ijkl} \nabla_{l} u.
\end{eqnarray}
Taking trace on $j$ and $k$, and using the second Bianchi identity
we have
\begin{eqnarray}
\nabla_{i} R -2 R_{ij} \nabla_{j} u =0,  \end{eqnarray} and
\begin{eqnarray}
\nabla_{i} (|\nabla u|^{2} + R -2\lambda u)=0. \end{eqnarray}

So, there is a constant $C_{9}$ such that
\begin{eqnarray}
|\nabla u|^{2} + R -2\lambda u = C_{9}. \end{eqnarray} As above,
we let $P\in M$ be the minimum point of $u$, then $|\nabla
u|(P)=0$, $\triangle u (P)\geq 0$, and $R(P)=n\lambda -\triangle u
(P)\leq n\lambda $. We have
\begin{eqnarray}
C_{9}=|\nabla u|^{2}(P) + R(P) -2\lambda u (P)\leq n\lambda
-2\lambda u(P).
\end{eqnarray}
From (2.12), we have
\begin{eqnarray}
|\nabla u|^{2}=- R +2\lambda u + C_{9}\leq n\lambda +2\lambda
(u-\inf_{x\in M}u)-R.
\end{eqnarray}

When the constant $\lambda \leq 0$, from the above inequality, we
have
\begin{eqnarray}
R\leq R+|\nabla u|^{2}\leq n\lambda ,
\end{eqnarray}
then $\triangle u =n\lambda -R \geq 0$. Since manifold $M$ is
compact, then $u$ must be a constant. So we have the following
Proposition.

\hspace{0.3cm}

{\bf Proposition 2.1 } {\it Let $g$ be a steady or expanding
gradient soliton over compact manifold $M$, then $g$ must be a
Einstein metric.}

\hspace{0.3cm}

When $\lambda$ is a positive constant, from the estimate (2.8) we
have
\begin{eqnarray}
|\nabla u|^{2}\leq n\lambda +2\lambda (\sup_{x\in M}u-\inf_{x\in
M}u)-R\leq C_{10},
\end{eqnarray}
where $C_{10}$ is a constant depending only on $C_{1}$, $C_{2}$,
$C_{3}$ and $\lambda$.

\hspace{0.3cm}

Let $(M_{\alpha}, g_{\alpha })$ be a sequence of shrinking Ricci
solitons  satisfying conditions (1), (2), (3) in theorem 1.1. From
(2.6), (2.8), (2.16), we obtain a uniform $C^{1}$-bound of
$u_{\alpha}$; from (2.8) and (2.14), we also obtain an uniform
upper bound of scalar curvature, i.e. we obtain the following
lemma,

\hspace{0.3cm}

{\bf Lemma 2.2 } {\it Let $(M_{\alpha}, g_{\alpha })$ be a
sequence of shrinking Ricci solitons ($\lambda =1$) satisfying
conditions (1), (2), (3) in theorem 1.1, and $u_{\alpha}$
satisfying the constraint (2.1), then there is positive constants
$C_{11}$, $C_{12}$ depending only on $C_{1}$, $C_{2}$ and $C_{3}$
such that
\begin{eqnarray}
|u_{\alpha}|_{C^{1}}\leq C_{11}.
\end{eqnarray}
and
\begin{eqnarray}
R(g_{\alpha})\leq C_{12}
\end{eqnarray}}

 \hspace{0.3cm}

In the next part of this section, we will give an uniform bound of
Perelman's functional $\mu (g , \frac{1}{2})$ ([Pe1]) for a
sequence shrinking Ricci solitons $(M_{\alpha}, g_{\alpha })$
satisfying conditions (1), (2), (3) in theorem 1.1. In [Pe1],
Perelman has introduced the following functional satisfying
\begin{eqnarray}
\textit{W}(g,\varphi , \tau )=(4\pi \tau )^{-\frac{n}{2}}\int_{M}
e^{-\varphi}[\tau (R+|\nabla \varphi |^{2})+f-n]\, dV_{g},
\end{eqnarray}
under the constraint
\begin{eqnarray}
(4\pi \tau )^{-\frac{n}{2}}\int_{M} e^{-\varphi}\, dV_{g}=1.
\end{eqnarray}
Then he define the functional
\begin{eqnarray}
\mu (g, \tau )=\inf \textit{W}(g, \cdot , \tau ),
\end{eqnarray}
where $\tau >0$, inf is taken over all functions satisfying the
constraint (2.20).

\hspace{0.3cm}

{\bf Lemma 2.3 } {\it If $(M, g)$ is a shrinking gradient Ricci
soliton, i.e. $$g-Ric(g)=\nabla du , $$ where $u$ satisfying the
constraint (2.1), then $u$ is a minimizer of Perelman's functional
$\textit{W}$ with respect to metric $g$ and $\tau =\frac{1}{2}$. }

\hspace{0.3cm}

{\bf Proof:} Let $\psi (t)$ be the $1$-parameter family of
diffeomorphisms that come from the vector field $\nabla u$, and
let $g(t)=\psi ^{\ast } g $, then $g(t)$ satisfy the following
Ricci flow equation,
\begin{eqnarray}
\frac{d}{dt}g(t)=-2Ric(g(t))+2g(t).
\end{eqnarray}
In order to use the Perelman's monotonicity formula [Pe1], we
scale the metric by $\tilde{g}(s)=C(s)g(t(s))$, where $C(s)=1-2s$,
$t(s)=-\frac{1}{2}\ln (1-2s)$. Then we have
\begin{eqnarray}
\frac{d}{ds}\tilde{g}(s)=-2Ric(\tilde{g}(s)),
\end{eqnarray}
and
\begin{eqnarray}
\mu (g(t), \frac{1}{2})=\mu (\tilde{g}(s(t)), \frac{1}{2}-s(t)).
\end{eqnarray}
 Let $\varphi (0)$ be a minimizer of $\textit{W}$ with respect
to metric $g(0)=g$ and $\tau =\frac{1}{2}$. Then function $\varphi
(t) =\psi ^{\ast }\varphi (0) $ is a minimizer of $\textit{W}$
with respect to metric $g(t)$ since
\begin{eqnarray}
\mu (g(t), \frac{1}{2})\leq \textit{W}(g(t), \varphi (t),
\frac{1}{2})=\textit{W}(g(0), \varphi (0), \frac{1}{2})=\mu (g(0),
\frac{1}{2})\leq \mu (g(t), \frac{1}{2}),
\end{eqnarray}
where the last inequality comes from the the Perelman's
monotonicity formula for $\mu (\tilde{g}(s), \frac{1}{2}-s)$.
Therefore, we have
\begin{eqnarray}
\begin{array}{lll}
0&=\frac{d}{dt}\textit{W}(g(t), \varphi (t), \frac{1}{2})\\
&=(2\pi )^{\frac{1}{2}}\int_{M} e^{-\varphi (t)}|R_{ij}+\varphi
_{ij}-g_{ij}|^{2}\, dV_{t},\\
\end{array}
\end{eqnarray}
which implies $\triangle \varphi (t)=n-R(t)=\triangle (u\circ \psi
(t) )$. Since $M$ is compact and both functions satisfy the the
constraint (2.1), then $\varphi (t)= u\circ \psi (t)$.

 $\Box$

\hspace{0.3cm}

For our sequence of shrinking gradient Ricci soliton $(M_{\alpha},
g_{\alpha} )$, the previous lemma tells us that every $u_{\alpha}$
is a minimizer of $\textit{W}(g_{\alpha }, \cdot , \frac{1}{2})$
and therefore satisfies ([Pe1])
\begin{eqnarray}
\triangle u_{\alpha} -\frac{1}{2}|\nabla u_{\alpha}|^{2}
+\frac{1}{2}R(g_{\alpha})+u_{\alpha}-n=\mu (g_{\alpha },
\frac{1}{2})
\end{eqnarray}
By $\triangle u_{\alpha}=n-R(g_{\alpha})$ and Lemma 2.2, we obtain
a uniform bound of $\mu (g_{\alpha}, \frac{1}{2})$.

\hspace{0.3cm}

{\bf Proposition 2.4} {\it Let $(M_{\alpha}, g_{\alpha })$ be a
sequence of shrinking Ricci solitons ($\lambda =1$) satisfying
conditions (1), (2), (3) in theorem 1.1, and $u_{\alpha}$
satisfying the constraint (2.1), then there is a constant $C_{13}$
depending only on $C_{1}$, $C_{2}$ and $C_{3}$ such that
\begin{eqnarray}
|\mu (g_{\alpha }, \frac{1}{2})|\leq C_{13}.
\end{eqnarray} }

\hspace{0.3cm}

{\bf Remark 2.5} {\it From Lemma 2.2 , we get scalar curvature
bounds for Ricci solitons satisfying conditions (1),(2),(3) in
theorem 1.1. Since Ricci curvature have a lower bound, we can also
get an upper bound for Ricci curvature.  If dimension $n=3$, we
know Ricci cuvature bounds imply Riemannian curvature bounds. On
the other hand, by Shi's estimates [Si] we can get  uniform bounds
for higher derivatives of Riemannian curvature, then using
Gromov-Cheeger compactness theorem, we can easily get the
following compactness theorem for Ricci solitons. }

\hspace{0.3cm}

{\bf Proposition 2.6} {\it Let $(M_{\alpha}, g_{\alpha})$ be a
sequence of shrinking gradient Ricci solitons of dimension $n=3$,
such that

(1),  $Ric(g_{\alpha})\geq -C_{1}g_{\alpha}$;

(2), $diam (M_{\alpha}, g_{\alpha})\leq C_{2}$;

(3), $Vol (M_{\alpha }, g_{\alpha })\geq C_{3}>0$;

for some uniform constants $C_{1}$, $C_{2}$, $C_{3}$. Then there
is a subsequence $(M_{\alpha}, g_{\alpha})$ converging to
$(M_{\infty}, g_{\infty})$ in $C^{\infty}$topology, and
$(M_{\infty}, g_{\infty})$ is a smooth gradient Ricci soliton.}

\hspace{0.3cm}

\section{$\epsilon $-regularity for Ricci solitons}
\setcounter{equation}{0}

\hspace{0.3cm}

Let $(M, g)$ be a shrinking gradient Ricci soliton. Choose a
normal coordinate system on the considered point, by direct
calculation, we have
\begin{eqnarray}
\begin{array}{lll}
&\triangle R_{ijkl} =\nabla _{m} \nabla_{m} R_{ijkl}\\
&=-\nabla _{m}\nabla _{k} R_{ijlm} -\nabla _{m}\nabla
_{l}R_{ijmk}\\
&=-\nabla _{k}\nabla _{m} R_{ijlm} -\nabla _{l}\nabla
_{m}R_{ijmk}+Q(Rm)_{ijkl}\\
&=\nabla _{k}\nabla _{m} R_{mlij} -\nabla _{m}\nabla
_{l}R_{mkij}+Q(Rm)_{ijkl}\\
&=\nabla _{k}\nabla_{i}R_{lj}-\nabla _{k}\nabla_{j}R_{li}-\nabla
_{l}\nabla_{i}R_{kj}+\nabla _{l}\nabla_{j}R_{ki}+Q(Rm)_{ijkl}\\
&=\nabla_{k}(R_{mlij}\nabla_{m}u)-\nabla_{l}(R_{mkij}\nabla_{m}u)+Q(Rm)_{ijkl}\\
&=\nabla_{k} R_{mlij}
\nabla_{m}u+R_{mlij}\nabla_{k}\nabla_{m}u-\nabla_{l}
R_{mkij}\nabla_{m}u-R_{mkij}\nabla_{l}\nabla_{m}u+Q(Rm)_{ijkl},\\
\end{array}
\end{eqnarray}
where we have used the second Bianchi identity, Ricci identity,
and formula (2.9), $Q(Rm)$ denotes a quadratic express in the
curvature tensor. In the short hand, we write the above identity
as
\begin{eqnarray}
\triangle Rm =\nabla Rm \ast \nabla u +Rm \ast g +Rm \ast Ric +Rm
\ast Rm.
\end{eqnarray}
Then, we have
\begin{eqnarray}
\begin{array}{lll}
\triangle |Rm|^{2}& =2|\nabla Rm|^{2} + 2 \langle \triangle Rm ,
Rm \rangle  \\
&\geq 2|\nabla Rm|^{2} -4|\nabla Rm||\nabla
u||Rm|-C_{14}|Rm|^{2}-C_{14}|Rm|^{3},\\
\end{array}
\end{eqnarray}
where $C_{14}$ is a positive constant depending only on dimension
$n$. By using the estimate (2.17) and Kato inequality, we get
\begin{eqnarray} \triangle |Rm|^{2}\geq (2-\theta ) |\nabla
|Rm||^{2}-C(\theta)|Rm|^{2}-C_{14}|Rm|^{3}.\end{eqnarray} Next, we
use the Moser iteration argument to deduce the following
mean-value inequality.

\hspace{0.3cm}

{\bf Lemma 3.1 } {\it  Let $(M, g)$ be a compact Riemannian
manifold, and $f$ be a Lipschitz function satisfying
\begin{eqnarray}
f\triangle f \geq -\theta_{1}|\nabla
f|^{2}-\theta_{2}f^{2}-\theta_{3}f^{3},
\end{eqnarray}
in weakly sense. Suppose that $\theta_{1}\leq\frac{1}{4}$, then
there exist constant $\epsilon $ depending only on the dimension
of $M$, $\theta_{3}$ and the lower bound of the Sobolev constant
$C_{s}$ so that  if
\begin{eqnarray}
\int_{B_{P}(2r)} f^{\frac{n}{2}} dv_{g}\leq \epsilon,
\end{eqnarray}
then
\begin{eqnarray}
\sup_{B_{P}(\frac{r}{2})}f\leq C_{\ast}
(1+\frac{1}{r^{2}})(\int_{B_{P}(r)}f^{\frac{n}{2}}\, dv_{g}
)^{\frac{2}{n}},\end{eqnarray} where $C_{\ast}$ depending only on
only on the dimension of $M$, $\theta_{2}$, $\theta_{3}$, the
lower bound of $Vol(M)$ and the Sobolev constant $C_{s}$ }

\hspace{0.3cm}

{\bf Proof:} Multiplying $\eta ^{2} f^{q-1}$ to (3.5), and
integrating yields
\begin{eqnarray}
\begin{array}{lll}
&\frac{4\theta_{1}}{q^{2}}\int_{M} \eta^{2}|\nabla
f^{\frac{q}{2}}|^{2}\, +\theta_{2} \int_{M} \eta^{2} f^{q}\,
+\theta_{3}\int_{M}\eta^{2}f^{q+1}\, \geq -\int_{M}\eta^{2}f^{q-1}\triangle f \,\\
&= \frac{4}{q}\int_{M}\eta f^{\frac{q}{2}} \langle \nabla \eta ,
\nabla f^{\frac{q}{2}} \rangle +\frac{4(q-1)}{q^{2}} \int_{M}
\eta^{2} |\nabla f^{\frac{q}{2}}|^{2}\,\\
&\geq -\frac{2}{q-1}\int_{M}f^{q}|\nabla \eta
|^{2}+\frac{2(q-1)}{q^{2}} \int_{M}
\eta^{2} |\nabla f^{\frac{q}{2}}|^{2}\, ,\\
\end{array}
\end{eqnarray}
where $q\geq 2$ and $\eta $ be a nonnegative cut off  function
that we will choose later. If we suppose that $\theta_{1}\leq
\frac{1}{4}$, from the above inequality, we have
\begin{eqnarray}
\frac{(q-1)}{q^{2}} \int_{M} \eta^{2} |\nabla
f^{\frac{q}{2}}|^{2}\,\leq \frac{2}{q-1}\int_{M}f^{q}|\nabla \eta
|^{2}+\theta_{2} \int_{M} \eta^{2} f^{q}\,
+\theta_{3}\int_{M}\eta^{2}f^{q+1}\, .
\end{eqnarray}
Using the Sobolev inequality (1.6), and let $\mu = \frac{n}{n-2}$,
we obtain
\begin{eqnarray}
\begin{array}{lll}
&\{\int_{M}  (\eta f^{\frac{q}{2}})^{2\mu}\, \}^{\frac{1}{\mu
}}\leq \frac{2}{C_{s}^{2}}\int_{M}|\nabla (\eta
f^{\frac{q}{2}})|^{2}\, +2Vol(M)^{-\frac{4}{n}}\int_{M}
\eta^{2}f^{q}\,\\
&\leq \frac{2}{C_{s}^{2}} \{\int_{M}\eta^{2} |\nabla
f^{\frac{q}{2}}|^{2}\, +\int_{M} f^{q}|\nabla \eta
|^{2}\}+2Vol(M)^{-\frac{4}{n}}\int_{M} \eta^{2}f^{q}\,\\
&\leq \frac{2\theta_{3}q^{2}}{C_{s}^{2}(q-1)}\int_{M}
\eta^{2}f^{q+1}\, +\frac{6q^{2}}{C_{s}^{2}(q-1)^{2}}\int_{M}
f^{q}|\nabla \eta
|^{2}\,\\
&+(\frac{2\theta_{2}q^{2}}{C_{s}^{2}(q-1)}+2Vol(M)^{-\frac{4}{n}})\int_{M} \eta^{2}f^{q}\, .\\
\end{array}
\end{eqnarray}
By H\"older inequality, we have
\begin{eqnarray}
\int_{M} \eta^{2}f^{q+1}\, \leq (\int_{Supp \eta
}f^{\frac{n}{2}}\,)^{\frac{2}{n}}\{\int_{M}(\eta
f^{\frac{q}{2}})^{2\mu }\,\}^{\frac{1}{\mu} }
\end{eqnarray}

 Take $\epsilon
<\{\frac{(n-2)C_{s}^{2}}{2\theta_{3}n^{2}}\}$. Let $q=\frac{n}{2}$
and $\eta $ be a cut off function with compact support in
$B_{P}(2r)$, equal to $1$ on $B_{P}(r)$ and such that $|\nabla
\eta |\leq \frac{2}{r}$, from the above inequalities, we get
\begin{eqnarray}
\{\int_{B_{P}(r)} f^{\frac{n}{2}\mu }\}^{\frac{1}{\mu }}\leq
(\frac{48n^{2}}{C_{s}^{2}(n-2)^{2}}\frac{1}{r^{2}}+\frac{4n^{2}}
{C_{s}^{2}(n-2)^{2}}\theta_{2}+4Vol(M)^{-\frac{4}{n}})\int_{B_{P}(2r)}
f^{\frac{n}{2}}\, .
\end{eqnarray}

In the next part of the proof, we choose cut off functions $\eta$
with compact support in $B_{P}(r)$, and equal to $1$ on
$B_{P}(\frac{r}{2})$. Using the H\"older inequality again, we get
\begin{eqnarray}
\int_{M} \eta ^{2} f^{q+1} \, \leq \{\int_{B_{P}(r)}
f^{\frac{n}{2}\mu }\}^{\frac{2}{n\mu }}(\int_{M}(\eta
f^{\frac{q}{2}})^{2\nu} )^{\frac{1}{\nu}}
\end{eqnarray}
where $\nu =\frac{\frac{n}{2}\mu }{\frac{n}{2}\mu -1
}=\frac{n^{2}}{n^{2}-2n+4}$. Using Young inequality, we have
\begin{eqnarray}
(\int_{M}(\eta f^{\frac{q}{2}})^{2\nu} )^{\frac{1}{\nu }}\leq
\delta (\int_{M}(\eta f^{\frac{q}{2}})^{2\mu}
)^{\frac{1}{\mu}}+C(n)\delta^{-\frac{n-2}{2}}\int_{M}\eta^{2}f^{q}\,
,
\end{eqnarray}
for small $\delta$.

Set $\delta = \frac{1}{2}\{\int_{B_{P}(r)} f^{\frac{n}{2}\mu
}\}^{\frac{2}{n\mu }}\frac{C_{s}^{2}(q-1)}{2\theta_{3}q^{2}}$,
from (3.10), (3.13) and (3.14), we have

\begin{eqnarray}
\{\int_{M} (\eta^{2} f^{q})^{\mu}\}^{\frac{1}{\mu }}\leq C_{15}
q^{n}\int_{M}((1+\frac{1}{r^{2}})\eta^{2}+|\nabla \eta
|^{2})f^{q}\, ,
\end{eqnarray}
where $C_{15}$ is a positive constant depending only on
$\theta_{2}$, $\theta_{3}$, dimension of $M$,  the lower bound of
$Vol(M)$ and the Sobolev constant.

Set $\frac{r}{2}\leq r_{2}< r_{1}\leq r $, and let $\eta \in
C_{0}^{\infty}(B_{P}(r_{1})) $ be the cut off function with the
property that $\eta =1$ in $B_{P}(r_{2})$ and $|\nabla \eta |\leq
\frac{2}{r_{1}-r_{2}} $. From (3.15), we have
\begin{eqnarray}
(\int_{B_{P}(r_{2})}  f^{q\mu})^{\frac{1}{\mu }}\leq 4C_{15}
q^{n}(1+\frac{1}{r^{2}}+\frac{1}{(r_{1}-r_{2})^{2}})\int_{B_{P}(r_{1})}f^{q}\,
.
\end{eqnarray}
Let $R_{i}=\frac{r}{2}+\frac{r}{2}2^{-i}$,
$q_{i}=\frac{n}{2}\mu^{i}$, applying (3.16) to $r_{1}=R_{i}$,
$r_{2}=R_{i+1}$, $q=q_{i}$, we have
\begin{eqnarray}
(\int_{B_{P}(R_{i+1})} f^{\frac{n}{2}\mu^{i+1}})^{\mu^{-(i+1)}
}\leq (64C_{15}
(1+\frac{1}{r^{2}})\frac{n}{2})^{\mu^{-i}}(2\mu)^{n\cdot
i\mu^{-i}}(\int_{R_{i}}f^{\frac{n}{2}\mu^{i}}\,)^{\mu^{-i}} .
\end{eqnarray}
Observe that $\lim_{i\rightarrow \infty }R_{i}=\frac{r}{2}$, and
iterating the above inequality, we conclude that
\begin{eqnarray}
\sup_{B_{P}(\frac{r}{2})}f^{\frac{n}{2}}\leq
C_{16}(1+\frac{1}{r^{2}})^{\frac{n}{2}}\int_{B_{p}(r)}f^{\frac{n}{2}}\,
.
\end{eqnarray}

$\Box$

\hspace{0.3cm}

Let $\theta =\frac{1}{2}$ in the formula (3.4), then the norm of
Riemannian curvature   $|Rm|$  of shrinking gradient Ricci
solitons must satisfy the Bochner type inequality (3.5). From
lemma 3.1, we obtain the $\epsilon$ regularity estimates for
shrinking gradient Ricci solitons.

\hspace{0.3cm}

{\bf Theorem 3.2 } {\it Let $(M_{\alpha}, g_{\alpha})$ be a
sequence of compact gradient Ricci solitons satisfying the
conditions (1), (2), (3) in the theorem 1.1. Then there exist
constants $C_{17}$ and  $\epsilon$ depending only on $C_{1}$,
$C_{2}$, $C_{3}$ so that if \begin{eqnarray}
\int_{B_{P}^{\alpha}(2r)}|Rm(g_{\alpha})|^{\frac{n}{2}}< \epsilon,
\end{eqnarray}
then
\begin{eqnarray}\sup_{B_{P}^{\alpha}(\frac{r}{2})}|Rm(g_{\alpha})|
\leq
C_{17}(1+\frac{1}{r^{2}})(\int_{B_{P}(r)}|Rm(g_{\alpha})|^{\frac{n}{2}}\,)^{\frac{2}{n}}.\end{eqnarray}}

\hspace{0.3cm}

\section{The proof of theorem 1.1}
\setcounter{equation}{0}

\hspace{0.3cm}

In this section, firstly, we will show that we can extract a
subsequence of Ricci solitons $(M_{\alpha}, g_{\alpha })$ which
satisfy conditions (1),(2),(3),(4) in theorem 1.1, so that it
converges to  an orbifold in a topological sense. This relies on
work by M.Anderson ([An1], [An2]).

Let $(M, g)$ be a Riemannian manifold, and $h_{M}$ be the
isoperimetric constant given by
\begin{eqnarray}
h_{M}=\inf_{S} \frac{(Vol(M))^{n}}{[\min (Vol(M_{1}),
Vol(M_{2}))]^{n-1}},
\end{eqnarray}
where $S$ varies over all closed hyper-surfaces of $M$ such that
$M\setminus S =M_{1}\cup M_{2}$. Croke ([Cr]) shows that $h_{M}$
is bounded below by a constant depending only on a lower bound for
Ricci curvature and  volume, and an upper bound on the diameter.
In particular, if $B_{x}(r)$ is a geodesic ball of radius $r$
about $x\in M$ and $S_{x}(r)=\partial B_{x}(r)$,
$v(r)=Vol(B_{x}(r))$, then it follows that: $(v' (r))^{n}\geq
h_{M}v(r)^{n-1}$, for $v(r)<\frac{1}{2}Vol(M)$, integrating this
inequality, one obtains: $v(r)\geq n^{-n}h_{M}r^{n}$. On the other
hand, from the Bishop volume comparison theorem, we know that
there must exist a positive constant $C_{18}$ depending only on
the lower bound of Ricci curvature and volume such that
\begin{eqnarray}
Vol(B_{x}(r))<\frac{1}{2}Vol(M), \quad  when  \quad r<C_{18}.
\end{eqnarray}
So, we have
\begin{eqnarray}
Vol(B_{x}(r))\geq C_{19} r^{n},
\end{eqnarray}
for $r<C_{18}$, where $C_{19}$ depend only on a lower bound for
the isoperimetric constant. From [Ya], volume noncollapsing
condition [4.3] and a lower bound for isoperimetric constant imply
the following Sobolev inequality
\begin{eqnarray}
(\int_{B_{x}(r)}f^{\frac{2n}{n-2}}\,)^{\frac{n-2}{n}}\leq
\frac{1}{C_{s}'}\int_{B_{x}(r)}|\nabla f|^{2}\, ,
\end{eqnarray}
for every Lipschitz function $f$ with compact support in
$B_{x}(r)$ and $r\leq C_{18}$. In fact, an upper bound on
diameter, a lower bound of Ricci curvature and volume
noncollapsing imply a lower bound on Sobolev constant $C_{s}'$.

\hspace{0.3cm}

Let $\epsilon$ be the constant of theorem 3.2 which determined by
the bounds $Ric(g_{\alpha})\geq -C_{1}g_{\alpha}$, $diam
(M_{\alpha }, g_{\alpha })\leq C_{2}$, $Vol(M_{\alpha},
g_{\alpha})\leq C_{3}$ on $(M_{\alpha}, g_{\alpha })$. We fix an
$0<r<C_{18}$ and let $\{x_{k}^{\alpha }\}$ be a maximal
$\frac{r}{8}$ separated set in $(M_{\alpha }, g_{\alpha})$. Thus
the geodesic balls $B_{x_{k}}^{\alpha}(\frac{r}{16})$ are
disjoint, and the balls $B_{x_{k}}^{\alpha}(\frac{r}{4})$ form a
cover of $M_{\alpha}$. We let
\begin{eqnarray}
G_{\alpha }^{r}=\cup \{B_{x_{k}}^{\alpha}(\frac{r}{4}):
\int_{B_{x_{k}}^{\alpha}(2r)}|Rm(g_{\alpha})|^{\frac{n}{2}}<\epsilon
\}, \nonumber
\end{eqnarray}
and
\begin{eqnarray}
B_{\alpha }^{r}=\cup \{B_{x_{k}}^{\alpha}(\frac{r}{4}):
\int_{B_{x_{k}}^{\alpha}(2r)}|Rm(g_{\alpha})|^{\frac{n}{2}}\geq
\epsilon \}. \nonumber
\end{eqnarray}
Then $M_{\alpha}=G_{\alpha }^{r}\cup B_{\alpha }^{r}$. From volume
non-collapsing condition (4.3) and Bishop-Gromov volume estimates,
we have a bound on the number of bad balls $Q_{\alpha }^{r}$ in
$B_{\alpha}^{r}$ independent of $\alpha $ and $r$, namely,
\begin{eqnarray}
Q_{\alpha}^{r}\leq C_{20}, \end{eqnarray} where $C_{20}$ is
positive constant depending only on the constants $C_{1}$,
$C_{2}$, $C_{3}$, $C_{4}$ which given in theorem 1.1.

Since Ricci solitons are the solutions of Ricci flow, Shi's
curvature estimates do apply and therefore by the estimates (3.20)
\begin{eqnarray}
\sup_{G_{\alpha }^{r}} |\nabla^{k} Rm(g_{\alpha})|\leq
\frac{C_{21}}{r^{k+2}},
\end{eqnarray}
where $C_{21}$ depend only on $C_{1}$, $C_{2}$, $C_{3}$. We also
obtain
\begin{eqnarray}
\sup_{G_{\alpha}^{r}}|\nabla^{k} u_{\alpha}|\leq C_{22}(k),
\end{eqnarray}
where $C_{22}(k)$ is a constant depending only on $k$, $r$,
$C_{1}$, $C_{2}$, $C_{3}$.

From the volume non-collapsing condition (4.3), small curvature
estimates (theorem 3.2 ), and Shi's curvature estimates (4.6),
following section 5 in [An1], we can show that there is a
subsequence of $(M_{\alpha}, g_{\alpha})$ converges to
$(M_{\infty}, g_{\infty})$ in the Hausdorff topology, and
$M_{\infty}=G\cup \{P_{i}\}_{1}^{Q} $ is a complete lenght space
with a length function $g_{\infty}$, which restricts to a smooth
gradient Ricci soliton on $G$ satisfying
\begin{eqnarray}
g_{\infty}-Ric(g_{\infty})=\nabla d u_{\infty},
\end{eqnarray}
where $u_{\infty}$ is a $C^{\infty}$ limit of $u_{\alpha}$ away
from singular points. The point $\{P_{i}\}_{1}^{Q}$ are called the
curvature singularities of $M_{\infty}$, and the convergence is in
$C^{\infty}$-topology outside the singularities. Then, in the
similar way as that in section 5 in [An1], we can check that
$M_{\infty}$ has the structure of an orbifold with a finite number
of curvature singularity points, each having a punctured
neighborhood which is diffeomorphic to a punctured cone on a
spherical space form, and metrics $g_{\infty}$ has a $C^{0}$
extension over every singularity points. So, we have proved the
following proposition.

\hspace{0.3cm}

{\bf Proposition 4.1 } {\it  Let $(M_{\alpha}, g_{\alpha})$ be a
sequence of compact gradient Ricci solitons satisfying the
conditions (1), (2), (3), (4) in the theorem 1.1. There is a
subsequence so that $(M_{\alpha}, g_{\alpha})$ converges to a
compact orbifold $(M_{\infty}, g_{\infty})$ with finitly many
singularities. Convergence is in $C^{\infty}$ topology outside
those singularity points, $g_{\infty}$ is a smooth Ricci soliton
outside those singularities and has a $C^{0}$ extension over every
singularity points.}

\hspace{0.3cm}

To finish the proof of Theorem 1.1 we still need to show that the
limit metric $g_{\infty}$ on $G$ can be extended to an orbifold
metric on $M_{\infty}$. More precisely, in an orbifold lifting
around singular points, in an appropriate gauge, the gradient
Ricci soliton equation of $g_{\infty}$ can be smoothly extended
over the origin in a ball in $R^{n}$. At this stage, the
regularity theory is not sufficient to imply that $g_{\infty}$ is
smooth. However, by Fatou's  lemma, we know that
\begin{eqnarray}
\int_{M_{\infty}}|Rm(g_{\infty})|^{\frac{n}{2}}\, dv_{\infty}
<\infty.
\end{eqnarray}
From the above inequality, we can obtain an upper bound for the
norm of curvature tensor $Rm(g_{\infty})$ of the limit metric
$g_{\infty}$.

\hspace{0.3cm}

{\bf Lemma 4.2 }  {\it $|Rm(g_{\infty})|_{\infty}$ is bounded
uniformly on $M_{\infty}\setminus \{P_{i}\}_{1}^{Q}$.}

\hspace{0.3cm}

We leave the proof of Lemma 4.2 in the next section. From above,
we know that each singular point $P\in M_{\infty}$ has a
neighbourhood that is covered by a punctured ball
$B^{n}(r)\setminus \{0\}\in R^{n}$. Our goal is to show that there
exist one diffeomorphism $\phi$ of $B^{n}(r)\setminus \{0\}$ such
that $\phi^{\ast}\pi^{\ast}(g_{\infty})$ extends to a smooth
metric on $B^{n}(r)$, where $\pi$ is the covering map.

Using Lemma 4.2, and harmonic coordinates constructed in [Jo], in
the same way as in [BKN, theorem 5.1], we can show that if $r$ is
sufficiently small, there is a diffeomorphism $\phi$ of
$B^{n}(r)\setminus \{0\}$ such that $\phi$ extends to a
homeomorphism of $B^{n}(r)$ and satisfies
\begin{eqnarray}
\begin{array}{lll}
(g_{\infty})_{ij}(x) -\delta _{ij} &=\textit{O}(|x|^{2}),\\
\partial_{k}(g_{\infty})_{ij}(x)&=\textit{O}(|x|),\\
\end{array}
\end{eqnarray}
where we also denote  the pulled back metric
$\phi^{\ast}\pi^{\ast}(g_{\infty})$ by $g_{\infty}$ for
simplicity. This means that there are some coordinates in a
covering of a singular point of $M_{\infty}$ in which $g_{\infty}$
extends to a $C^{1,1}$-metric.

For our sequence of shrinking gradient Ricci soliton $(M_{\alpha},
g_{\alpha} )$, the lemma 2.3 tells us that every potential
function $u_{\alpha}$ is a minimizer of $\textit{W}(g_{\alpha },
\cdot , \frac{1}{2})$ and therefore satisfies
\begin{eqnarray}
\triangle u_{\alpha} -\frac{1}{2}|\nabla u_{\alpha}|^{2}
+\frac{1}{2}R(g_{\alpha})+u_{\alpha}-n=\mu (g_{\alpha },
\frac{1}{2}). \nonumber\end{eqnarray} By $\triangle
u_{\alpha}=n-R(g_{\alpha})$, we have \begin{eqnarray} \triangle
u_{\alpha }=|\nabla u_{\alpha}|^{2}-2u_{\alpha}+n+2\mu
(g_{\alpha}, \frac{1}{2}).
\end{eqnarray}
Proposition 2.4 tell us that those $\mu (g_{\alpha },
\frac{1}{2})$ are bounded uniformly. So we can extract a
subsequence of a sequence of converging metrics $g_{\alpha}$ such
that
\begin{eqnarray}
\lim_{\alpha \rightarrow \infty }\mu (g_{\alpha },
\frac{1}{2})=\mu_{\infty}.
\end{eqnarray}
Let $\alpha \rightarrow \infty$ in (4.11) we get
\begin{eqnarray} \triangle
u_{\infty }=|\nabla u_{\infty}|^{2}-2u_{\infty}+n+2\mu_{\infty},
\end{eqnarray}
away from those singular points $P_{i}$. On the other hand, from
Lemma 2.2 and the soliton equation , we have
\begin{eqnarray}
\sup_{M_{\infty}\setminus \{P_{i}
\}_{1}^{Q}}|u_{\infty}|_{C^{2}}\leq C_{\ast},
\end{eqnarray}
for some uniform constant $C_{\ast}$. So, it is not hard to
conclude that $\nabla u_{\infty}$ extends to the origin in the
covering ball $B^{n}(r)$. Moreover, $u_{\infty}\in
C^{1,1}(B^{n}(r))$.

Using the harmonic coordinates for $g_{\infty}$ in $B^{n}(r)$, we
can write the soliton equation as follow
\begin{eqnarray}\triangle (g_{ij})+ \cdots =u_{ij}\end{eqnarray}
where the dots indicate lower order terms involving at most one
derivative of $g_{ij}$. Since $g_{\infty}\in C^{1,1}(B^{n}(r))$
and $u_{\infty}\in C^{1,1}(B^{n}(r))$, from (4.11) and (4.15), the
standard elliptic regularity theory imply that $g_{\infty}$ and
$u_{\infty}$ must be smooth in $B^{n}(r)$, that is $g_{\infty}$ is
a gradient soliton metric in an orbifold sense.

When $n$ is odd, we can discuss like that in  [An1, section5], to
conclude that there are no curvature singularities in
$M_{\infty}$. We argue by contradiction. Suppose that there exist
curvature singularities in $M_{\infty}$.  For each curvature
singularity $P\in \{P_{i}\}_{1}^{Q}\subset M_{\infty}$, there is a
sequence $x_{\alpha}\in M_{\alpha}$, such that
$x_{\alpha}\rightarrow P$ and $\inf_{r>0}\sup \{|Rm_{\alpha}(x)|;
x\in B_{x_{\alpha}}(r)\subset M_{\alpha}\}\rightarrow \infty$, as
$\alpha \rightarrow \infty$. Since the curvature of $M_{\alpha}$
remains bounded in bounded distance away from $x_{\alpha}$, we may
assume that $x_{\alpha}$ realizes the maximum $R_{\alpha}$ of
$|Rm_{\alpha}|$ on $B_{x_{\alpha}}(r_{0})$ for some small $r_{0}$.
Now consider the pointed connected Riemannian manifolds
$V_{\alpha}=(B_{x_{\alpha}}, x_{\alpha},
R_{\alpha}^{\frac{1}{2}}g_{\alpha})$. We note that the curvature
of $V_{\alpha}$ is uniformly bounded,
$|Rm_{V_{\alpha}}(x_{\alpha})|=1$, and
$|Ric(V_{\alpha})|(x)\rightarrow 0$ for any point $x\in
V_{\alpha}$ as $\alpha \rightarrow \infty$ (since,in Lemma2.2 we
have proved that the Ricci curvature of $M_{\alpha}$ is bouunded
uniformly ).  Similarly, $\int_{V_{\alpha }}|Rm|^{\frac{n}{2}}\leq
C$ and the Sobolev constants for $V_{\alpha}$ are uniformly
bounded below, since this is true for $M_{\alpha}$ itself. As in
section three in [An1], we can prove that there is  a subsequence
of $V_{\alpha}$ converges, in $C^{\infty}$ topology on compact
sets, to a complete connected Riemannian manifold $V$ satisfying
\begin{eqnarray}
\begin{array}{lll}
Ric_{V}&=0,\\
\frac{Vol(B(r))}{r^{n}}&\geq C',\\
\int_{V}|Rm|^{\frac{n}{2}}\, &\leq C,
\end{array}
\end{eqnarray}
and
\begin{eqnarray}
|Rm|(x_{0})=1, \quad for some x_{0}\in V.
\end{eqnarray}
A complete connected Riemannian manifold satisfying (4.16) is
called an EALE (Einstein, asymptotically locally Euclidean) space.
Theorem 3.5 in [An1] had shown that in odd dimensions, nontrivial,
i.e., nonflat, EALE space do not exist, so we get the
contradiction by (4.17). So it follows that $M_{\infty}$ is a
smooth manifold and the convergence $M_{\alpha}\rightarrow
M_{\infty}$ is smooth.

\hspace{0.3cm}

\section{Curvature bounds of the limiting metric}
\setcounter{equation}{0}

\hspace{0.3cm}

In this section, we will give curvature bounds for the limiting
metric $g_{\infty}$. If $n=4$, the approach that we will use to
prove lemma 4.2 is based on Uhlenbeck's [Uh, Th4.1] idea for
treating the isolated singularities for the Yang-Mills equation.
If $n\geq 5$, using the Sobolev bounds in $B_{g_{\infty}}(P, r)\
\{P\}$, we can verify the basic methods of Sibner [Sib, Lemma2.1,
Proposition2.4] remain valid here also. We should point out that
Chao and Sesum [CS] had used the same idea to treating
the¡¡K\"ahler-Ricci soliton case. We include a sketched proof
here.

Let $P$ be a singular point of $M_{\infty}$,
$r(x)=dist_{g_{\infty}}(x , P)$, and $B_{g_{\infty}}(P, r)=\{x\in
M_{\infty}| r(x)<r \}$. Since the Sobolev inequality (4.4) with a
uniform Sobolev sonstant $C_{s}'$ for all $g_{\alpha}$, we have
the following Sobolev inequality
\begin{eqnarray}
(\int_{B_{g_{\infty}}(P,
r)}f^{\frac{2n}{n-2}}\,)^{\frac{n-2}{n}}\leq
\frac{1}{C_{s}'}\int_{B_{g_{\infty}}(P, r)}|\nabla f|^{2}\, ,
\end{eqnarray}
for every Lipschitz function $f$ with compact support in
$B_{g_{\infty}}(P, r)\setminus \{P \}$ and $r\leq C_{18}$.

\hspace{0.3cm}

{\bf Remark 5.1 } {\it By Fatou's lemma and similar arguments as
in [BKN], we can get that (5.1) also holds for compact supported
functions $f\in W^{1,2}(B_{g_{\infty}}(P, r))$.}

\hspace{0.3cm}

 From
(3.4), we also have
\begin{eqnarray} \triangle |Rm|^{2}\geq (2-\theta ) |\nabla
|Rm||^{2}-C(\theta)|Rm|^{2}-C_{14}|Rm|^{3},\end{eqnarray} on
$M_{\infty}\setminus \{P_{i}\}_{1}^{Q}$. by (4.9), we can decrease
$r$ if necessary so that $\int_{g_{\infty}(P,
r)}|Rm|^{\frac{n}{2}}\, dV_{g_{\infty}}<\epsilon$, where $\epsilon
$ is chosen to be small. By Sobolev inequality (5.1) and lemma
3.1, we have
\begin{eqnarray}
|Rm(g_{\infty})|(x)\leq
\frac{C}{r(x)^{2}}\{\int_{B_{g_{\infty}}(P, 2r(x)
)}|Rm(g_{\infty})|^{\frac{n}{2}}\,
dV_{g_{\infty}}\}^{\frac{2}{n}},
\end{eqnarray}
for some uniform constant $C$.  From Shi's curvature estimates
(4.6), letting $\alpha \rightarrow \infty$ we get
\begin{eqnarray}
|\nabla^{k}Rm(g_{\infty})|(x)\leq \frac{C(r(x))}{r(x)^{k+2}},
\end{eqnarray}
for all $x\in M_{\infty}\setminus \{P_{i}\}_{1}^{Q}$, where
$C(r(x))\rightarrow 0$ as $r(x)\rightarrow 0$.

 \hspace{0.3cm}

(a), When $n=4$. Let $U$ be a small neighbourhood of $P$, recall
that $U\setminus \{P\}$ is covered by $B^{n}(r)\setminus
\{0\}\subset R^4$ and $\pi^{\ast}g_{\infty}$ extends to a $C^{0}$
metric on the ball $B^{n}(r)$, where $\pi$ is the covering map. By
estimates (5.3) and (5.4), as in [Ti1, section 4] we can find a
gauge $\phi $ (i.e, a diffeomorphism on $B^{n}(r)(r)$) so that
\begin{eqnarray}
\begin{array}{ll}
|d g_{ij}|(x) &\leq \frac{\epsilon (r(x))}{r(x)},\\
|d (\frac{\partial g_{ij}}{\partial x_{k}})|(x) &\leq \frac{\epsilon (r(x))}{r(x)^{2}},\\
|d (\frac{\partial ^{2} g_{ij}}{\partial x_{k}\partial x_{l}})|(x) &\leq \frac{\epsilon (r(x))}{r(x)^{3}},\\
\end{array}
\end{eqnarray}
in $B^{r}\setminus \{0\}$, where $d$ is the exterior differential
on $R^{4}$ and $|\cdot |$ is the norm with respect to the
Euclidean metric, and $g$ stands for
$\phi^{\ast}\pi^{\ast}g_{\infty}$.

Let $D=d+A$ be a connection  uniquely associated to the metric $g$
on $B^{n}(r)\setminus \{0\}$, where $A$ is the connection form. As
in [Uh, section4] or [Ti1, section4], we constructed the broken
Hodge gauges. We break domain up into annuli
\begin{eqnarray}
\begin{array}{lll}
\textbf{U}_{l}&=\{ x : 2^{-l-1}r\leq r(x)\leq 2^{-l}r\},\\
\textbf{S}_{l}&=\{ x : r(x)= 2^{-l}r\},
\end{array}
\end{eqnarray}
for $l=0,1,2,\cdots$.

\hspace{0.3cm}

{\bf Definition 5.2} {\it A broken Hodge gauge for a connection
$D$ in a bundle $E$ over $\cup_{l=0}^{\infty}\textbf{U}_{l}$ is a
gauge related continuously to the original gauge in which $D=d+A$
and $A(l)=A|_{\textbf{U}_{l}}$ have the following properties for
all $l$ \begin{eqnarray} \begin{array}{lll} &d^{\ast}A(l)=0 \quad
in \textbf{U}_{l},\\&
A_{\psi}(l)|_{\textbf{S}_{l}}=A_{\psi}(l-1)|_{S_{l}},\\
&d_{\psi}^{\ast} A_{\psi} (l)=0 \quad on \textbf{S}_{l} \quad and
\textbf{S}_{l+1},\\&\int_{\textbf{S}_{l}}A_{r}(l)\,
=\int_{\textbf{S}_{l+1}}A_{r}(l)=0.\\
\end{array}\nonumber\end{eqnarray} }

\hspace{0.3cm}

As in [Uh, Theorem4.6] or [Ti1,section4], from estimate (5.3), we
have the following lemma

\hspace{0.3cm}

{\bf Lemma 5.3} {\it Let $D$ be the unique connection associate to
the metric $g$, for small $r$,then there exists a broken Hodge
gauge in $B^{n}(r)\setminus
\{0\}=\cup_{l=0}^{\infty}\textbf{U}_{l}$  satisfying
\begin{eqnarray}\begin{array}{lll}& |A(l)|_{g}(x)\leq C 2^{-l}r\sup_{\textbf{U}_{l}}|Rm|_{g}\leq C
2^{l+1}r^{-1},\\&\int_{\textbf{U}_{l}} |A(l)|_{g}^{2}\, dV_{g}\leq
C 2^{-2l}r^{2}\int_{\textbf{U}_{l}}|Rm|_{g}^{2}\,
dV_{g}.\\\end{array}\end{eqnarray}}

\hspace{0.3cm}

By direct calculation, we have
\begin{eqnarray}
\begin{array}{lll}
\sum_{l=0}^{\infty}\int_{\textbf{U}_{l}}|Rm|_{g}^{2}\,
dV_{g}&=-\sum_{l=0}^{\infty}\int_{\textbf{U}_{l}} \langle A(l),
D^{\ast}Rm \rangle\,
-\sum_{l=0}^{\infty}\int_{\textbf{U}_{l}} \langle [A(l), A(l)], Rm\rangle\,\\
&+\int_{\textbf{S}_{0}}\langle A_{\psi}(0),
Rm_{r\psi}\rangle\,-\lim_{l\rightarrow
\infty}\int_{\textbf{S}_{l+1}}\langle A_{\psi}(l),
Rm_{r\psi}\rangle\,\\
\end{array}
\end{eqnarray}

From (5.3) and (5.7), it is not hard to check that
$\lim_{l\rightarrow \infty}\int_{\textbf{S}_{l+1}}\langle
A_{\psi}(l), Rm_{r\psi}\rangle\,=0$. On the other hand, we know
the limit metric $g_{\infty}$ satisfies the Ricci soliton
equation,
\begin{eqnarray}
g-Ric=\nabla du,
\end{eqnarray}
where $u=\phi^{\ast}\pi^{\ast}u_{\infty}$. We have
\begin{eqnarray}
D^{\ast}Rm_{ijk}=R_{ijkm, m}=R_{ki,
j}-R_{kj,i}=u_{k,ji}-u_{k,ij}=R_{ijkl}g^{lm}u_{m}.
\end{eqnarray}
where we have used the second Bianchi identity and Ricci identity.
By lemma 2.2, we know that $|\nabla u|$ is bounded uniformly, so
we have
\begin{eqnarray}
\int_{\textbf{U}_{l}} \langle A(l), D^{\ast}Rm \rangle\,\leq
(\int_{\textbf{U}_{l}}
|A(l)|^{2}\,)^{\frac{1}{2}}(\int_{\textbf{U}_{l}}
|D^{\ast}Rm|^{2}\,)^{\frac{1}{2}}\leq
C2^{-l}r\int_{\textbf{U}_{l}}|Rm|_{g}^{2}\, .
\end{eqnarray}
From the estimate (5.3), we get
\begin{eqnarray}
\begin{array}{lll}
&|\int_{\textbf{U}_{l}} \langle [A(l), A(l)], Rm\rangle\,|\leq
\sup_{\textbf{U}_{l}}|Rm|_{g}\int_{\textbf{U}_{l}} |A(l)|^{2}\,\\
&\leq C(\int_{B_{g_{\infty}}(P,
2^{-l+1}r)}|Rm(g_{\infty})|^{2}\,)^{\frac{1}{2}}\int_{\textbf{U}_{l}}
|Rm|_{g}^{2}\,\\.
\end{array}
\end{eqnarray}
Similar as the proof of corollary 2.6 in [Uh] with small
modification just like that in [Ti1, setion4], we can find a
decreasing function $\epsilon_{1} (r)$ with $\lim_{r\rightarrow
0}\epsilon_{1} (r)=0$ such that
\begin{eqnarray}
\int_{\textbf{S}_{0}}|A_{\psi}(0)|^{2}_{g}\, dV_{g}\leq
(2-\epsilon (r))^{-2}r^{2}\int_{\textbf{S}_{0}}|Rm_{\psi \psi
}|^{2}_{g} dV_{g}
\end{eqnarray}

From (5.8), (5.11), (5.12), (5.13), we have
\begin{eqnarray}
\begin{array}{lll}
\int_{B(r)}|Rm|_{g}^{2} \, dV_{g}&\leq
\frac{r}{2(2-\epsilon_{1}(r))(1-\epsilon_{2}(r))}\int_{\partial
B(r)}|Rm|_{g}^{2} \,dV_{g}\\
&\leq \frac{r}{4}(1+\frac{\delta}{2})\int_{\partial
B(r)}|Rm|_{g}^{2}\, dV_{g},\\
\end{array}
\end{eqnarray}
whenever $r$ is sufficiently small and $\delta \in (0, 1)$, where
$\lim_{r\rightarrow 0}\epsilon_{2}(r)=0$ . Then it is standard to
conclude from above inequality that([Ti1, section4])
\begin{eqnarray}
|Rm|_{g_{\infty}}(x)\leq \frac{C}{r(x)^{\delta}},
\end{eqnarray}
for $x\in B_{g_{\infty}}(P, r)$, for sufficiently small $r$ and
some $\delta \in (0, 1)$, where $C$ is a uniform constant. Recall
that $g_{\infty}$ extends $C^{0}$ metric on the covering ball,
from (5.15), we can show that there exist $q>4$ so that
\begin{eqnarray}
\int_{B_{g_{\infty}}(P, r)}|Rm(g_{\infty})|^{q}\, <\infty .
\end{eqnarray}

To further consideration, we need the following lemma which
similar as lemma 2.1 in [Sib], and had been proved in [CS]

\hspace{0.3cm}

{\bf Lemma 5.4} {\it Let $f\geq 0$ be a smooth function in
$M_{\infty}\setminus \{P_{i}\}_{1}^{Q}$ and satisfies (3.5), with
$f\in L^{\frac{n}{2}}$. If $f\in L^{\frac{2nq_{0}}{n-2}}\cap
L^{2q}$, $q_{0}>\frac{1}{2}$, then $\nabla f^{q}\in L^{2}$ and for
sufficiently small $r$, we have
\begin{eqnarray}
\int_{B_{g_{\infty}}(P_{i}, r)}\eta^{2}|\nabla f^{q}|^{2}\, \leq
\int_{B_{g_{\infty}}(P_{i}, r)} |\nabla \eta |^{2} f^{2q},
\end{eqnarray}
for all $\eta \in B_{g_{\infty}}(P_{i}, r)$, where $C$ is a
uniform constant.

}

\hspace{0.3cm}

(b), If $n>4$, let $f=|Rm(g_{\infty})|\in L^{\frac{n}{2}}$, we can
choose $q_{0}=1$ and $q=\frac{n}{4}$. From (5.2), we know
$f=|Rm(g_{\infty})|\in L^{\frac{n}{2}}$ satisfies (3.5), applying
Lemma 5.4 to $f$, we find that $\nabla f^{\frac{n}{4}}\in L^{2}$.
By remark 5.1, we can apply the Sobolev inequality (5.1) to
$f^{\frac{n}{4}}$ to conclude that
\begin{eqnarray}
|Rm(g_{\infty})|\in L^{p},
\end{eqnarray}
with $p=\frac{n}{2}\frac{n}{n-2}>\frac{n}{2}$.

\hspace{0.3cm}

From above we know that $|Rm(g_{\infty})|\in L^{p}$ for some
$p>\frac{n}{2}$ in both cases $n=4$ and $n>4$. Specially, since
$Vol(M_{\infty}, g_{\infty})<\infty$, by (5.16) (5.18) and using
H\"older inequality, we have $|Rm(g_{\infty})|\in L^{p}$ for $p\in
(0, \frac{n}{2}\frac{n}{n-2}]$. Take $q_{0}=1$, $q\in (0,
\frac{n}{4}\frac{n}{n-2}]$ and repeat applying Lemma 5.4 to get
$\nabla |Rm_{g_{\infty}}|^{q} \in L^{2} $. By Sobolev inequality
(5.1) and Using H\"older inequality again, we have
$|Rm(g_{\infty})|\in L^{p}$ for $p\in (0,
\frac{n}{2}(\frac{n}{n-2})^{2}]$. If we keep on repeating this, at
the $k$-th step we have $\nabla |Rm(g_{\infty})|^{q}\in L^{2}$ for
$q\in (0, \frac{n}{4}(\frac{n}{n-2})^{k}]$ and
$|Rm(g_{\infty})|\in L^{p} $ for $p\in (0,
\frac{n}{2}(\frac{n}{n-2})^{k+1}]$. Since
$(\frac{n}{n-2})^{k}\rightarrow \infty$ as $k\rightarrow \infty$,
we have \begin{eqnarray}|Rm(g_{\infty})|\in L^{p}, \end{eqnarray}
and \begin{eqnarray} \nabla |Rm(g_{\infty})|^{p} \in L^{2}, \quad
for \quad all\quad p.
\end{eqnarray}

By (5.19) (5.20), combining Remark 5.1 and lemma 5.4, for
sufficiently small $r$ and any $p$ we have
\begin{eqnarray}
\begin{array}{lll}
(\int_{B_{g_{\infty}}(P, r)}|\eta
|Rm_{g_{\infty}}|^{p}|^{\frac{2n}{n-2}}\,)^{\frac{n-2}{n}}&\leq
\frac{1}{C_{s}'}\int_{B_{g_{\infty}}(P, r)} |\nabla (\eta
|Rm(g_{\infty})|^{p})|^{2}\,\\&\leq C\int_{B_{g_{\infty}}(P,
r)}|\nabla \eta |^{2}|Rm(g_{\infty})|^{2p}\, ,\\
\end{array}
\end{eqnarray}
with a uniform constant $C$, where $\eta $ is any cut off function
with compact support in $B_{g_{\infty}}(P, r)$. Then, using the
Moser's iteration argument as in the proof of Lemma 3.1, we get
\begin{eqnarray}
\sup_{B_{g_{\infty}}(P, \frac{r}{2})}|Rm(g_{\infty})|\leq
\frac{C}{r^{2}}.
\end{eqnarray}
So, we get the curvature bound for the limiting metric
$g_{\infty}$and finished the proof of Lemma 4.2.

\section{References}
\setcounter{equation}{0}

\hspace{0.4cm}

\bigskip

{\parindent=0pt

\def\toto#1#2{\centerline{\hbox to 1.6cm{#1\hss}
\parbox[t]{158mm}{#2}}\vspace{\baselineskip}}

\toto{[An1]}{M.Anderson: Ricci curvature bounds an Einstein
metrics on compact manifolds, {\it J.Amer.Math.Soc.} {\bf 3}(1990)
355-374.}\toto{[An2]}{M.Anderson: Orbifold compactness for spaces
of Riemannian metrics and applications, {\it Math.Ann.} {\bf
331}(2005), 739-778. }\toto{[Be]}{A.Besse, Einstein manifolds,
{\it Ergeb.Math.Grenzgeb.Band 10, Springer-Verlag, Berlin and New
York}, 1987.}\toto{[BKN]}{S.Bando, A.Kasue, H.Nakajima: On a
construction of coordinates at infinity on manifolds with fast
curvature decay and maximal volume growth. {\it Invent.Math.} {\bf
97 }(1989), 313-349. }\toto{[Ch]}{J.Cheeger: Finiteness theorem
for Riemannian manifols. {\it Am. J. Math.} {\bf 92}(1970),
61-74.} \toto{[CS]}{H.D.Cao and N.Sesum: The compactness result
for K\"ahler Ricci solitons;  {\it
arXiv.math.DG/0504526}.}\toto{[Ga]}{Z.Gao: Einstein metrics, {\it
J.Differential Geometry} {\bf 32}(1990)
155-183.}\toto{[Gr]}{M.Gromov: Structures metriques pour les
varietes Riem. {\it Red. par J.Lafontaine et P.Pansu ,} Paris
(19881).}\toto{[GW]}{R.Greene and H.Wu: Lipschitz convergence of
Riemannian manifolds. {\it Pac.J.Math.} {\bf 131}(1988),
119-141.}\toto{[Ha1]}{R.Hamilton: Three-manifolds with positive
curvature Ricci curvature,{\it J. Differential Geometry} {\bf
17}(1982) 225-306.}\toto{[Ha2]}{R.Hamilton: Four-manifolds with
positive curvature operator, {\it J.Diff.Geom.} {\bf 24}(1986),
153-179. }\toto{[Jo]}{J.Jost: Harmonic mappings between Riemannian
manifolds. {Proceedings of the center for mathematical analysis},
Australian National University, Vol.{\bf 4}, 1983.
}\toto{[Ka]}{A.Kasue: A convergence theorem for Riemannian
manifolds and some applications. {\it Nagoya Math. J.} {\bf
114}(1989), 21-51.}\toto{[Na]}{H.Nakajima: Hausdorff convergence
of Einstein 4-manifolds, {\it J.Fac.Sci.Univ.Tokyo} {\bf 35}(1988)
411-424.}\toto{[Pe1]}{G.Perelman: The entropy formula for the
Ricci flow and its geometric application, {\it
arXiv.math.DG/0211159 }. }\toto{[Pe2]}{G.Perelman: Perelman's
lectures on the K\"ahler Ricci flow.}\toto{[Pe]}{S.Peters:
Convergence of Riemannian manifolds. {\it Compos.Math.} {\bf
62}(1987), 3-16.}\toto{[Si]}{W.X.Shi: Deforming the metric on
complete Riemannian manifolds; {\it J.Differential Geometry}, {\bf
30 } (1989) 223-301. }\toto{[Sib]}{L.M.Sibner: Isolated point
singularity problem; {\it Math.Ann.} {\bf 271}(1985) 125-131.
}\toto{[Ti1]}{G.Tian: On Calabi's conjecture for complex surface
with positive first Chern class; {\it Inventiones Math. } {\bf
101} (1990), 101-172.} \toto{[Ti2]}{G.Tian: Compactness theorems
for Kahler-Einstein manifolds of dimension $3$ and up, {\it
J.Differential Geometry} {\bf 35} (1992) 535-558.}
\toto{[Uh]}{K.Uhlenbeck: Removable singularities in Yang-Mills
field, {\it Comm. Math. Phys.} {\bf 83}(1982)
11-29.}\toto{[Ya]}{S.T.Yau: Survey lecture, {\it Seminar on
Differential Geom.}, {\it Ann.of Math. Stud.}, {\bf 102}(1982).}
\toto{[Zh]}{X.Zhang: Harnack inequality and regularity of
p-Laplace equation on complete manifolds, {\it Kodai Math.J.} {\bf
23 }(2000), 326-344.}
\bigskip

\end{document}